\documentclass[a4paper, 12t]{article}
\usepackage[latin1]{inputenc}
\usepackage[english]{babel}

\usepackage{amsmath}
\usepackage{amsthm}
\usepackage{amssymb}
\usepackage{color}
\usepackage{bbm}
\usepackage{wasysym}

\usepackage{tikz,pgfplots,subcaption}

\theoremstyle{plain}
\newtheorem{theorem}{Theorem}

\newtheorem{remark}[theorem]{Remark}

\newcommand{\R}{\mathbb{R}}

\newcommand{\mL}{\mathcal{L}}

\newcommand{\dd}{\, \text{d}}

\usepackage{hyperref}
\usepackage{vmargin}
\setmarginsrb{3cm}{2.5cm}{3cm}{1.5cm}{0cm}{0cm}{0cm}{1.5cm}

\newcommand{\matlab}{MATLAB\textsuperscript{\textregistered}}

\newlength\fheight
\newlength\fwidth
\begin{document}

\title{Neural Network Based Nonlinear Observers}

\maketitle

\centerline{\scshape Tobias Breiten}
\medskip
{\footnotesize
 \centerline{Institute of Mathematics and Scientific Computing}
   \centerline{University of Graz}
   \centerline{Graz, Austria}
} 

\medskip

\centerline{\scshape Karl Kunisch}
\medskip
{\footnotesize
 \centerline{Institute of Mathematics and Scientific Computing}
   \centerline{University of Graz}
   \centerline{Graz, Austria}
   \centerline{RICAM Institute}
   \centerline{Austrian Academy of Sciences}
   \centerline{Linz, Austria}
}

\bigskip

\begin{abstract}
  Nonlinear observers based on the well-known concept of minimum energy estimation are discussed. The approach relies on an output injection operator determined by a Hamilton-Jacobi-Bellman equation and is subsequently approximated by a neural network. A suitable optimization problem allowing to learn the network parameters is proposed and numerically investigated for linear and nonlinear oscillators.
\end{abstract}

\section{Introduction}
This paper is concerned with the problem of estimating the current state  $x(T)\in \R^n$ of a nonlinear finite dimensional system, given past observations  $y(s)\in \R^r$, for $s\in (0,T), T>0$, with $r$ typically smaller than $n$, and some a-priori knowledge $x_0$ of the initial condition. More precisely we are considering a pair consisting of a state and an observation process of the form
\begin{equation}\label{eq:system}
\begin{aligned}
\dot{x} & = f(x) + Gv, \ x(0)=x_0+\zeta, \\
 y & = Cx + w,
\end{aligned}
\end{equation}
where $v$ represents unobserved disturbance in the dynamics, $w$ stands for unobserved disturbance in the observation and $\zeta$ is an unobserved disturbance of the initial condition. The model dynamics $f$, and the matrices  $G$ and $C$ are assumed  to be known.
In the linear case, with $f(x)=Ax,\; A\in \R^{n\times n}$ the theory of dynamic observers is highly developed. Roughly speaking one can rely on either Luenberger observers, which build a dynamical system of the form
\begin{equation}\label{eq:Luen}
\begin{aligned}
\dot{\hat x} = A \hat x + L (y- C\hat x), \ \hat{x}(0)=x_0, \\
\end{aligned}
\end{equation}
for an appropriately chosen $L$ guaranteeing qualitative properties are guaranteed. Most notably one aims for convergence of the
observation error $x(t)-\hat x(t)$  to $0$ as $\to  \infty$. Alternatively one can focus  on  optimal observers which achieve the estimation of the unknown state while simultaneously minimizing a conveniently chosen variational principle. In the linear case this  leads to the celebrated Kalman-Bucy filter and the construction of $L$ in the observer equation \eqref{eq:Luen} involves computing the solution of a matrix Riccati equation.

In the nonlinear case, the construction of observers is significantly more complex. It appears to be the case that there is no first choice which is universally practiced. For the analogue of Luenberger observers the properties of $f$ will significantly influence the structure and construction of the analogue of $L$ in \eqref{eq:Luen} above, which needs to be done in a case by case approach. The Kalman analogue almost inevitable leads to the need of solving some kind of a Hamilton-Jacobi-Bellman (HJB) equation. The practical realization of nonlinear optimal observer theory or state reconstruction is therefore impeded by the curse of dimensionality.

 In our work we aim at making a step forward in developing optimal observers for nonlinear systems. For this purpose, we rely on the concept of the Mortensen observer \cite{Mortensen68}, also known as \emph{minimum energy estimation}. The concept is quite involved and therefore we considered it useful to recall the key steps of the derivation within this paper. We shall also provide a novel addition  by arguing that, contrary to its first appearance which suggests the observer gain  (corresponding to $L$ above) depends on $y$, this is in fact not the case. This brings the gain of the Mortensen observer closer to the Kalman filter gain than was suggested by our understanding of the literature.
 The characterization of this observer gain, however, still involves the inverse of the Hessian of the solution of a HJB equation. Thus computing the gain on the basis of the HJB equation is a challenging task, especially if one ultimately aims at investigating  systems of sizable dimension.
 To alleviate some of the difficulties we propose to approximate the gains by a family of parameterized functions. These parameters will be learned (optimized) by means of solving an ensemble of optimization problems which  approximate the Mortensen minimum energy functional. We have some freedom in the choice of the parameterized functions which approximate the observer gains. Here we choose to use a neural network  to profit from their good approximation properties.

The  theory of mathematical observers for linear systems is well-covered in many monographs. Here we can only list a small selection \cite{BucJ68,Jaz70,KalB61,KK85,Son98}.
We also mention the inspiring survey article \cite{Willems04}. Concerning contributions from nonlinear observer theory we mostly rely on \cite{Fleming97,Krener04,Moireau18,Mortensen68}. The use of neural networks for nonlinear observers is not new. They are used, however, for an entirely different purpose than in our work, namely for the approximation of the HJB equation itself, see e.g. \cite{Ad12}. To the contrary we bypass to directly solve the HJB equation, but we take advance of the relationship between the derivative of the value function and the construction of the observer gain.

Let us briefly outline the structure of the paper. In section 2 the relevant concepts from nonlinear observer theory are summarized. We do not claim any originality, but  this section can be of interest, since the material was not so straightforwardly available from other sources. Section 3 sets forth the learning approach that we propose to train the neural networks which approximate the observer gains. The last section contains first numerical results which demonstrate, first, to which degree the proposed approach coincides with Kalman filtering results in the case of linear problems, and, second, the feasibility and success for a nonlinear problem.

In this paper we aim at presenting our ideas for constructing observers based on training neural networks. For the most part the mathematical development is formal. It is our aim to provide detailed analysis and rigorous proofs, as well as to advance code development for treating higher dimensional problems,  in independent future research steps.

\section{Basic concepts}
\subsection{The Mortensen observer}
To describe the Mortensen observer let 
us commence  by considering
\begin{equation}\label{eq:P_T}\tag{$P^T$}
\begin{cases}
\begin{array}{rl}
&\min\limits_{\substack{v\in L^2(0,T;\R^m) \\ w \in L^2(0,T;\R^r)}}
  ~J(T,\xi; v, w) = \frac{1}{2} \|x(0) - x_0\|^2_{Q_0} + \frac{1}{2} \int_{0}^{T} (\|v(s)\|^2 + \alpha \|w(s)\|^2) ~\mathrm{d}s \\[2ex]
&\dot{x} = f(x) + Gv, ~x(T) = \xi,\\
&y = Cx + w,
\end{array}
\end{cases}
\end{equation}
where $Q_0=Q_0^\top \succ 0$ symmetric positive definite  and $y\in L^2(0,T;\R^r)$ are given. Further $G\in \R^{n\times m},\, C\in \R^{r\times n}$,  $f\colon\R^n\to \R^n$ is a $C^2-$ mapping, and $\|x\|_{Q_0}^2= x^\top Q_0 x$. Throughout we assume that the indicated minima exist.
\eqref{eq:P_T} is the basic building block for the Mortensen-based nonlinear filtering approach. In \cite{Mortensen68} the functional $J$ in \eqref{eq:P_T} is referred to as likelihood functional while in \cite{Krener04} it is called the minimum energy functional.

We  set
\begin{equation}\label{eq:val}
V(T,\xi) = \underset{(v,w)}{\min} J(T,\xi; v, w).
\end{equation}
and refer to
\begin{equation}\label{eq:observer}
\hat{x}(T) = \underset{\xi \in \mathbb{R}^n}{\arg \min} ~V(T,\xi)
\end{equation}
as the  minimum norm (Mortensen) observer estimate of the state at time $T$.
In \eqref{eq:val} we assumed that \eqref{eq:observer} admits a unique solution.
Let $x^T(\cdot)=x^T(\cdot;\hat x(T),v^T)$, denote the trajectory satisfying the first equation in \eqref{eq:system} with terminal condition $\hat x (T)$ and associated optimal $v^T$.

In the following developments of this subsection we draw from  \cite{Fleming97, Krener04, Moireau18, Mortensen68, Willems04}. It has been observed and will be detailed in Section \ref{sec:kalman} below that, in the linear case, the Mortensen observer  is closely related to the Kalman filter based state reconstruction, see e.g. \cite[Lemma 1, Theorem 3] {Willems04}.

For the further development and in particular to derive a characterizing equation for $\hat x(T)$ it is necessary to also consider for  all $t \in [0,T]$ the problems
\begin{equation}\label{eq:P_t} \tag{$P^t$}
\begin{cases}
\begin{array}{rl}
&\underset{v\in L^2(0,t;\R^m)}{\min} ~J(t,\xi; v) = \frac{1}{2} \|x(0) - x_0 \|^2_{Q_0} + \frac{1}{2} \int_{0}^{t} (\|v(s)\|^2 + \alpha \|y(s)-Cx(s)\|^2) ~\mathrm{d}s \\[2.5ex]
&\dot{x} = f(x) + Gv, ~x(t) = \xi, \\
\end{array}
\end{cases}
\end{equation}
with the associated value function given by
\begin{equation*}
V(t,\xi) = \underset{v}{\min}\, J(t,\xi; v).
\end{equation*}
Note that \eqref{eq:P_t} for $t=T$ coincides with \eqref{eq:P_T} where the observation  process is eliminated from the constraints and directly realized within the cost functional. Analogous to \eqref{eq:observer} we introduce
\begin{equation}\label{eq:observert}
\hat{x}(t) = \underset{\xi \in \mathbb{R}^n}{\arg \min} ~V(t,\xi),
\end{equation}
and the associated optimal trajectories $x^t(\cdot)=x^t(\cdot;\hat x(t),v^t)$ on $[0,t]$.
Let us observe that $\hat x(t)= x^t(t;\hat x(t),v^t)$ for every $t\in [0,T]$, but
$\hat x$ restricted to $[0,t]$, does not coincide with $ x^t(\cdot;\hat x(t),v^t)$ on $[0,t]$.  We next aim at deriving an equation for the mapping $t\mapsto \hat x(t)$.

Applying Bellman's principle to \eqref{eq:P_t} we obtain
\begin{equation*}
V(T,\xi)= V(t,x^*(t)) +  \underset{v\in L^2(t,T;\R^m)}{\min}\frac{1}{2} \int^T_t  (\|v(s)\|^2 + \alpha \|y(s)-Cx(s)\|^2) ~\mathrm{d}s +\frac{1}{2}\|x^*(t) - x(t)\|^2_{Q_0} ,
\end{equation*}
where $x$ is the solution to
\begin{equation}\label{eq:aux1}
\dot{x} = f(x) + Gv,  \text{ on } (t,T),~x(T) = \xi,
\end{equation}
and $x^*$ is the solution to the minimization problem on $[t,T]$. Note that compared to the familiar optimal control related Bellman principle \eqref{eq:aux1} has the format of being \emph{backwards in time}.

Assuming $C^1$-regularity of $V$, it satisfies the following HJB equation
\begin{equation*} \label{HJB} \tag{HJB}
\begin{aligned}
\begin{cases}
& \partial_t {V}(t,\xi) = - \nabla_{\xi} V(t, \xi)^\top f(\xi) - \frac{1}{2} \|G^\top \nabla_{\xi}V(t,\xi)\|^2 + \frac{\alpha}{2} \|y(t) - C\xi\|^2,\\[1.5ex]
& V(0,\xi) = \frac{1}{2} \|\xi - x_0\|_{Q_0}^2.
\end{cases}
\end{aligned}
\end{equation*}
The verification theorem implies that an optimal control $v$ in feedback form for \eqref{eq:P_t} is given by
\begin{equation}\label{eq:veri}
v(s)=  G^\top\nabla_\xi V(s,x(s)), \text{ for } s \in (0,t),
\end{equation}
see e.g.\cite[Theorem I.5.1.]{FS06}.

To present a governing equation for the observer $\hat x$ we first note that $T$ in \eqref{eq:observer}  was arbitrary and thus we can use \eqref{eq:observer} as the characterizing property for $\hat x(t)$ for every $t\in [0,T]$. This implies that
\begin{equation}\label{eq:character}
0 = \nabla_\xi V(t,\hat{x}(t)) \text{ for all } t\in [0,T].
\end{equation}
From \eqref{eq:character} we deduce that
\begin{equation}\label{eq:property}
\frac{\dd}{\dd t} V(t,\hat{x}(t)) =   \partial_t V(t,\hat{x}(t)) \; \text{ and }\; \partial_t \nabla_\xi V(t,\hat x(t))+ \nabla_{\xi \xi} V(t,\hat x(t)) \dot{\hat{ x}}(t)=0,  \text{ for all } t\in [0,T],
\end{equation}
where $\nabla_{\xi \xi} V(t,\cdot)$ denotes the Hessian of $V(t,\cdot)$. From the first equations in \eqref{HJB} and \eqref{eq:property} we obtain
\begin{equation}\label{eq:derivval}
\frac{\dd}{\dd t} V(t,\hat{x}(t)) = \frac{\alpha}{2}\|y(t)-C\hat x(t)\|^2,
\end{equation}
as observed in \cite{Krener04, Mortensen68}.  Taking the spatial derivative in the first equation of \eqref{HJB} along $(t,\hat x(t))$ we obtain
\begin{equation*}
0= \partial_t \nabla_\xi V(t,\hat x(t)) + \nabla_{\xi\xi}V(t,\hat x(t))f(\hat x(t)) + \alpha C^\top(y(t)-C \hat x(t)),
\end{equation*}
and combined with the second equation in \eqref{eq:property}
\begin{equation*}
0= - \nabla_{\xi \xi} V(t,\hat x(t)) \dot{\hat{ x}}(t) + \nabla_{\xi\xi}V(t,\hat x(t))f(\hat x(t)) + \alpha C^\top(y(t)-C \hat x(t)).
\end{equation*}

Assuming the invertibility of $\nabla_{\xi \xi} V$ along $(t,\hat x(t))$ and denoting
$\Pi(t)= (\nabla_{\xi \xi} V(t, \hat{x}(t)))^{-1}$
we arrive at the observer equation:

\begin{equation}\label{eq:obseq}
\begin{cases}
\begin{array}{l}
\dot{\hat{x}} = f(\hat{x}(t)) + \alpha \Pi(t)C^\top(y(t) - C \hat{x}(t)),\\[1.3ex]
\hat x(0)=x_0,
\end{array}
\end{cases}
\end{equation}
where the initial condition is a consequence of \eqref{eq:observert} for $t=0$.

Next we aim for deriving an equation for $\Pi$. For this purpose we take two derivatives of \eqref{HJB} with respect to the spatial variable and evaluate along $(t,\hat x(t))$ to obtain:
\begin{equation*}
\begin{array}l
 \partial_t \nabla_{\xi\xi}V(t,\hat x(t)) + Df(\hat x(t))^\top \nabla_{\xi \xi} V(t,\hat x(t)) +  \nabla_{\xi \xi} V(t,\hat x(t)) Df(\hat x(t)) \\[1.8ex]
+ \sum_{k=1}^n \partial_{\xi_i}\partial_{\xi_j}\partial_{\xi_k} V(t,\hat x(t)) f_k(\hat x(t)) +  \nabla_{\xi \xi} V(t,\hat x(t))G G^\top \nabla_{\xi \xi}V(t,\hat x(t)) - \alpha C^\top C=0,
 \end{array}
\end{equation*}
where $Df$ denotes the Jacobian of $f$.


Multiplying this equation by $\Pi(t)$ from the left and the right, we obtain,
\begin{equation}\label{eq:Piequation}
\partial_t{\Pi}(t)\!=\! Df(\hat{x}(t))\Pi (t) + \Pi(t) Df(\hat{x}(t))^\top\! +\! G G^\top\! -\! \alpha \Pi C^\top C \Pi\! +\! \Pi \sum_{k=1}^n \partial_{\xi_i}\partial_{\xi_j}\partial_{\xi_k} V(t,\hat x(t)) f_k(\hat x(t)) \Pi,
\end{equation}
where  $\partial_t{\Pi}(t)$ stands for $\partial_t (\partial_{\xi \xi}V(t,\hat x(t)))^{-1}$.  While we do not further use \eqref{eq:Piequation} it will be convenient to compare it to the Kalman filter equation further below.




\subsection{Relationship to Lagrangian formulation}
Associated to \eqref{eq:P_t} we introduce the Lagrange functional
\begin{align*}
\mL&\colon W^{1,2}(0,t;\R^n) \times L^{2}(0,t;\R^m)\times L^2(0,t;\R^n)\to \R , \\
\mL(x,v,p) &= \frac{1}{2} \|x(0) - x_0 \|^2_{Q_0} + \frac{1}{2} \int_{0}^{t} (\|v(s)\|^2 + \alpha \|y(s)-Cx(s)\|^2) ~\dd s \\
&\qquad + \int_0^t p(s)^\top (f(x(s)) + Gv(s)-\dot x(s))\, \dd s.
\end{align*}
It is straightforward to argue that the Lagrange multiplier rule is applicable. Hence we can derive the following first order optimality condition for \eqref{eq:P_t}:
\begin{equation}\label{eq:OC}
 \begin{aligned}
   \dot{x}_t(s) &= f(x_t(s)) + Gv_t(s),  \quad x_t(t)=\xi, \\[1.8ex]
     -\dot{p}_t(s)  &= Df(x_t(s))^\top p_t(s) -\alpha C^\top(y(s)-Cx_t(s)),\quad p_t(0)= - Q_0(x_t(0)-x_0),\\[1.8ex]
v_t(s)&=- G^\top p_t(s),
 \end{aligned}
\end{equation}
where $s\in (0,t)$. Note that $x_t$ does not coincide with $x^t$, which is the optimal trajectory for the minimum norm terminal $\xi$ according to \eqref{eq:observert}.  We also have the following relationship between the gradient of the value function and the adjoint variable $p$:
\begin{equation}\label{dualrelation}
p_t(s)= -\nabla_\xi V (s,x_t(s)) \text{ for every } t\in (0,T] \text{ and } s\in (0,t],
\end{equation}
see e.g.\cite[Theorem I.6.2]{FS06}.

\subsection{The Kalman-Bucy filter}\label{sec:kalman}

In the linear case, i.e., $f(x)=Ax$, we can give an explicit expression for $V(t,\xi)$ by means of a differential Riccati equation, see, e.g., \cite{Moireau18}.  For the purpose of a self-contained presentation, we detail the required  calculations. Let us define for $t\ge 0$ and $\xi \in \R^n$
\begin{align}\label{eq:V_linear}
 V(t,\xi) &:= \frac{1}{2} (\xi-\hat{x})^\top  \Sigma^{-1} (\xi-\hat{x}) + \frac{\alpha}{2} \int_0^t \| y-C\hat{x}\|^2\; \dd s,
\end{align}
where
\begin{align*}
 \dot{\Sigma} &= A\Sigma+\Sigma A^\top  - \alpha \Sigma C^\top  C \Sigma + GG^\top  , \ \Sigma(0)= Q_0^{-1}, \\
 \dot{\hat{x}}&= A\hat{x} + \alpha \Sigma C^\top  (y-C\hat{x}), \ \hat{x}(0)=x_0.
\end{align*}
In the above formulas the temporal dependence of the variables $\hat x$, $y$, and $\Sigma$ was suppressed. The Riccati equation satisfied by $\Sigma$ is well-known, see \cite{KalB61}. Let us note that in the current case, the Riccati equation and \eqref{eq:Piequation} coincide since the term $\sum_{k=1}^n \partial_{\xi_i}\partial_{\xi_j}\partial_{\xi_k} V(t,\hat x(t)) $ vanishes by the ansatz \eqref{eq:V_linear}. We further note the following properties
\begin{align*}
\frac{\mathrm{d}}{\mathrm{d}t} \Sigma^{-1} &= -\Sigma^{-1} \left(\frac{\mathrm{d}}{\mathrm{d}t} \Sigma \right) \Sigma^{-1} =  -\Sigma^{-1} A - A^\top \Sigma^{-1} + \alpha C^\top  C - \Sigma^{-1} GG^\top \Sigma^{-1}, \\
\nabla _\xi V(t,\xi) &=\Sigma^{-1} (\xi-\hat{x}), \ \ \nabla_{\xi}^2 V(t,\xi)= \Sigma^{-1}, \\
\| y-C\xi \|^2 &= \| (y-C\hat{x})-C(\xi-\hat{x})\|^2 = \| y-C\hat{x}\|^2 - 2  ( y-C\hat{x})^\top C(\xi-\hat{x}) + \| C(\xi-\hat{x}) \|^2.
\end{align*}
Subsequent use of the latter relations, allows us to show that $V$ solves the HJB equation:
\begin{align*}
\partial_t V(t,\xi) &= -\dot{\hat{x}}^\top \Sigma^{-1} (\xi -\hat{x})+\frac{1}{2}(\xi -\hat{x})^\top
\dot{\Sigma}^{-1} (\xi -\hat{x})  +\frac{\alpha}{2} \| y-C\hat{x}\|^2 \\
&= -(A\hat{x}+\alpha \Sigma C^\top (y-C\hat{x}))^\top \Sigma^{-1} (\xi -\hat{x}) + \frac{\alpha}{2} \| y-C\hat{x}\|^2  \\
& \quad -\frac{1}{2} (\xi-\hat{x})^\top (\Sigma^{-1} A + A^\top \Sigma^{-1}-\alpha C^\top C + \Sigma^{-1} GG^\top \Sigma^{-1} ) (\xi -\hat{x}) \\
&= -\alpha (y-C\hat{x})^\top C (\xi -\hat{x}) + \frac{\alpha}{2} \| y-C\hat{x} \|^2 +\frac{\alpha}{2} \|C (\xi -\hat{x})\|^2 \\
&\quad -\hat{x}^\top A^\top \Sigma^{-1} (\xi -\hat{x}) - \frac{1}{2} (\xi-\hat{x})^\top (\Sigma^{-1} A + A^\top \Sigma^{-1} + \Sigma^{-1} GG^\top \Sigma^{-1} ) (\xi -\hat{x}) \\
&=\frac{\alpha}{2} \| y-C\xi \|^2  -\frac{1}{2} \| G^\top \nabla _{\xi} V(t,\xi) \|^2 \\&\qquad  -\hat{x}^\top A^\top \nabla_\xi V(t,\xi)-\frac{1}{2}\nabla_\xi V(t,\xi)^\top A (\xi-\hat{x})-\frac{1}{2} (\xi-\hat{x})^\top A^\top \nabla_\xi V(t,\xi)  \\
&= \frac{\alpha}{2} \| y-C\xi \|^2  -\frac{1}{2} \| G^\top \nabla _{\xi} V(t,\xi) \|^2 -(A\xi )^\top \nabla_\xi V(t,\xi).
 \end{align*}
This is the first equation in \eqref{HJB}.  Further $V(0,\xi) := \frac{1}{2} (\xi-x_0)^\top  Q_0(\xi-x_0)$  as desired. In view of \eqref{eq:obseq}   the operator $\Sigma$ resumes the role  of $\Pi(t)= (\nabla_{\xi \xi} V(t, \hat{x}(t)))^{-1}$ in the linear case. The superposition  $\Sigma C^\top$ is is well-known from Kalman-Bucy filter.

 \subsection{$y$-independence of the observer gain in the nonlinear case}
 \label{subsec:y-indp}

 In the linear case, we can use the explicit expression for $V(t,\xi)$ to argue that $\nabla_\xi V(t,\xi)$ as well as $\nabla_\xi^2 V(t,\xi)$ are independent of $y$. In particular, this leads to an appealing feature of the optimal Kalman filter: varying the observation $y$ (by e.g., varying noise in initial state, dynamics, or output) does not require to recompute the observer gain $\Sigma C^\top$.
 In the nonlinear case this is not clear a-priori, since  $\Pi$ is constructed from $V$, which itself depends on $y$ as described in \eqref{HJB}. 
 In the following, we argue that in  the nonlinear case as well, the observer gain is independent of $y$. With regard to the learning approach in the subsequent section, this will allow us to design a ``global'' observer gain that is optimal for multiple observations $y$ at once.

 We recall  the HJB equation
  \begin{align*}
  \begin{cases}
  \partial_ t V(t,\xi)&=-\nabla_\xi V(t,\xi)^\top f(\xi) - \frac{1}{2} \| G^\top \nabla_\xi V(t,\xi) \|^2 + \frac{\alpha}{2} \| y(t) - C\xi \|^2, \\[1.5ex]
  V(0,\xi) &= \frac{1}{2} \| \xi - x_0 \|_{Q_0}^2.
 \end{cases}
 \end{align*}
Suggested by the definition of $V$ in the linear case,  let us define $W(t,\xi)$ by
\begin{align*}
 W(t,\xi):= V(t,\xi ) - \frac{\alpha}{2} \int_0^t \| y(s) - C x^{*,t}(s)\|^2 \, \mathrm{d}s,
\end{align*}
where $x^{*,t}$ is the solution to \eqref{eq:P_t} with $\xi\in \R^n$, it thus satisfies
\begin{align*}
\begin{cases}
 \dot{x}^{*,t}(s)&= f(x^{*,t}(s)) + G G^\top \nabla_\xi V(s, x^{*,t}(s)) \; \text{ for } s\in (0,t],\\[1.5ex]
 x^{*,t}(t)&=\xi.
 \end{cases}
\end{align*}
We have the equalities
\begin{align}\label{eq:gradVW}
 \nabla_\xi W(t,\xi) = \nabla _\xi V(t,\xi) , \ \  \nabla^2_\xi W(t,\xi) = \nabla^2 _\xi V(t,\xi).
\end{align}
With the previous relations, we now obtain
\begin{align*}
 \partial_t W(t,\xi) &= \partial_t V(t,\xi) - \frac{\alpha}{2} \| y(t) -C x^{*,t}(t)\|^2 \\
 &= -(\nabla_\xi W(t,\xi))^\top f(\xi) - \frac{1}{2} \| G^\top \nabla_\xi W(t,\xi) \|^2 + \frac{\alpha}{2} \| y(t) - C\xi \|^2 - \frac{\alpha}{2} \| y(t) - Cx^{*,t}(t) \| ^2 \\
 &= -(\nabla _\xi W(t,\xi))^\top f(\xi) - \frac{1}{2} \| G^\top \nabla \xi W(t,\xi) \|^2, \quad W(0,\xi)=V(0,\xi).
\end{align*}
Thus $W(t,\xi)$ is independent of $y$, and by \eqref{eq:gradVW} so are $\nabla_\xi V(t,\xi)$ and $\nabla _\xi ^2 V(t,\xi)$. In the subsequent section, we propose a learning approach that aims at approximating $\nabla_\xi V(t,\xi)$ and $\nabla _\xi ^2 V(t,\xi)$. The above considerations in particular allow us to design an observer that can be used for varying observations $y$ without re-learning the observer gains.

\section{Neural network based approximations of the Mortensen observer}

\subsection{A learning formulation for observer design}

From \eqref{eq:obseq} of the previous section it is evident that the observer gain of the Mortensen observer depends on the inverse of the Hessian of the solution to an HJB equation. The practical realization is therefore a formidable task. This motivates us to propose an approach which builds on the structural properties of the Mortensen observer but which does not depend on the availability of the  solution to the HJB equation \eqref{HJB}. For short, we shall replace  the disturbance in the dynamics $v$  by a neural network. This network will be trained by using \eqref{eq:P_t}, and information on the value function will be recovered on the basis of \eqref{eq:veri}. A related network based approach for stabilization of nonlinear systems was recently proposed in \cite{KW20}.

To describe the neural network formulation that we  propose in more detail, let us recall two of the main structural equations  of the previous section: following \eqref{eq:veri} and \eqref{eq:obseq} the closed-loop optimal solution to \eqref{eq:P_T}, and the observer equation are given by
\begin{equation}\label{eq:aux2}
\begin{cases}
\begin{array}{ll}
\dot{x}(t) = f(x(t)) + GG^\top\nabla_\xi V(t,x(t)),  ~x(T) =\xi, \\[1.7ex]
\dot{\hat{x}}(t) = f(\hat{x}(t)) + \alpha \Pi(t)(y(t) - C \hat{x}(t)), ~\hat{x}(0) = x_0, \text{ for } t \in (0,T).
\end{array}
\end{cases}
\end{equation}
These equations  depend on the gradient, respectively the inverse Hessian of the value  function $V(t,\xi)$ related to \eqref{eq:P_t}. In the following we describe a methodology which approximates these mappings by a network based function.

More specifically we approximate $\nabla_\xi V(t,x)$ as

\begin{equation}\label{eq:feedappr}
\nabla_\xi V(t,x)\approx h_\theta(t,x) = g_\theta(t,x) - g_s(t),
\end{equation}
where $g_\theta\colon[0,T]\times \R^n \to \R^n$, denotes a parameter dependent family of functions. In our case $\theta$ denotes the network parameters and $g_\theta$ will have the structure of a neural network. The shift $g_s\in L^2(0,T;\R^n)$ is a user defined function, which is chosen such that it reflects \eqref{eq:character}. We could choose $g_s(t)=g_{\theta}(t, \hat{x}(t))$, with the shift function $g_s$ depending on the network parameter. This would imply that $h_\theta(t,\hat{x}(t))=0$, independently of $\theta$. However,  the optimal estimate $\hat x$ is unknown a-priori. Moreover, in the course of characterizing and implementing optimality conditions with respect to $\theta$ further below, this would require to compute additional derivatives. Therefore, in \eqref{eq:feedappr} we propose to use a  function depending on time only, to introduce a shift to the network $g_\theta$. Note that this shift is not an additional approximation, it is simply our choice of making an ansatz for the approximation to $\nabla_\xi V(t,x)\approx h_\theta(t,x)$.
We will readdress this topic in the numerical section and discuss potential strategies for choosing $g_s$.

These considerations lead  to approximations of \eqref{eq:aux2} given by
\begin{equation}
\begin{cases}
\begin{array}{rl}
&\dot{x}_\theta(t) = f(x_\theta(t)) + G G^\top h_\theta(t,x_\theta(t)), \; x(T)=\xi, \\
&\dot{\hat{x}}_\theta (t)= f(\hat{x}_\theta(t)) + (D_x h_\theta(t, \hat{x}_\theta(t)))^{-1} C^\top(y(t) - C \hat{x}_\theta(t)), \; \hat{x}_\theta(0) = x_0.
\end{array}
\end{cases}
\end{equation}

Here the degrees of freedom represented by $\theta$  should be constructed such that $h_\theta$ is a good approximation of $\nabla_\xi V(t,x)$, which in turn provides an approximation  for $v$ in \eqref{eq:P_T} in feedback form. Recall here that $w$ in  \eqref{eq:P_T} can equivalently be expressed as $y-Cx$.

The network parameters $\theta$ will be determined by considering \eqref{eq:P_T} with $v$ replaced by $G^\top h_\theta$. If this was the only information for determining $\theta$ it could suffer from the fact that it would depend too strongly on a particular choice of the terminal state $\xi$. For this reason we  choose an ensemble  $\{\xi_j\}_{j=1}^d$ of possible terminal states  in $\R^n$. The parameters $\theta$ are then determined by solving:
\begin{equation}\label{eq:Ptheta}\tag{$P_{T,\theta}$}
\left\{
\begin{aligned}
 \min_{\substack {\theta \in \mathbb R^N\\ x_{\theta,j} \in W^{1,2}(0,t;\mathbb R^n)} } J(\theta,x_{\theta})&:= \frac{1}{d}\sum\limits_{j=1}^d \big (  \, \frac{1}{2} \|x_{\theta,j}(0) - x_0\|_{Q_0}^2 + \frac{1}{2} \int\limits_{0}^{T} (\|G^\top h_\theta(t, x_{\theta,j}(t))\|^2 + \alpha \|y(t)- C x_{\theta,j}(t)\|^2) \, \dd t  \big)\\[1.8ex]
\text{s.t.\ }\dot{x}_{\theta,j}(t) &= f(x_{\theta,j}(t))+ G G^\top h_\theta(t,x_{\theta,j}(t)),\;  x_{\theta,j}(T) = \xi_j, \ j=1,\dots,d.
\end{aligned}\right.
\end{equation}

\subsection{Parametrization by neural networks}

In this subsection, we collect some notation that is standard in the context of neural networks, see e.g. \cite{GKNV19,GPEB19}.  For the construction of a function $g_\theta\colon[0,T]\times \mathbb R^n \to \mathbb R^n$ let us fix $L\in \mathbb N$ as well as  $n_i\in \mathbb N$ for $ i=0,\dots,L$ and consider a parameter set $\theta$ given by
\begin{align*}
 \theta = (\theta_1,\dots,\theta_L)=(W_1,b_1,R_1,\dots,W_{L-1},b_{L-1},R_{L-1},W_L),
\end{align*}
with $W_i,R_i \in \mathbb R^{n_i\times n_{i-1}}$ and $b_i \in \mathbb R^{n_i}$. We then define $g_\theta$ as a composition of functions as follows
\begin{align*}
&  g_\theta \colon \mathbb R^{n+1} \to \mathbb R^{n}, g_{\theta} (z)=\left(g_{\theta_L}\circ g_{\theta_{L-1}} \circ \cdots \circ g_{\theta_1} \right)(z), \\
 & g_ {\theta_i}\colon \mathbb R^{n_{i-1}} \to \mathbb R^{n_i}, \
  g_{\theta_i}(z)=\sigma (W_i z + b_i)+R_i z, \ \ i=1,\dots, L-1, \\
  &g_ {\theta_L}\colon \mathbb R^{n_{L-1}} \to \mathbb R^{n_L}, \
  g_{\theta_L}(z)=W_L z ,
\end{align*}
where $n_0=n+1$ and $n_L=n$. Here, the evaluation of the \textit{activation function} $\sigma$ is defined componentwise  via
\begin{align*}
 \sigma\colon \mathbb R^\ell \to \mathbb R^{\ell},  (\sigma(z))_i=\sigma(z_i).
\end{align*}
In our examples below, we choose $\sigma$ as the \emph{logistic function}, i.e., $\sigma(s)=\frac{1}{1+\exp(-s)}$. However, plenty of other choices are conceivable.  The matrices $R_i$ are sometimes referred to as \textit{residual connection}   and their use has been shown to be beneficial in particular with respect to numerical stability, see, e.g., \cite{HeZRS16}.

\subsection{Learning the parameters via optimization}

The learning process for the function $g_{\theta}$ then is characterized by the optimization problem \eqref{eq:Ptheta} for $(\theta,x_{\theta,j}) \in  \mathbb R^N \times W^{1,2}(0,t;\mathbb R^n)$, where $N=n_L\cdot n_{L-1}+\sum_{i=1}^{L-1} (2 \cdot n_{i-1}+1) n_{i} $. In this respect we can refer to a rather detailed description in \cite{KW20} on the treatment of the optimization problem which describes the network learning step.
Assuming the existence of a minimizer  $(\theta^*,x_{\theta,j}^*),j=1,\dots,d$ and associated Lagrange multipliers $p_{\theta,j},j=1,\dots,d$, we can formally  characterize $(\theta,x_{\theta_j},p_{\theta,j})$ as solutions to the first-order optimality conditions
\begin{equation}\label{eq:network_1st_opt}
\begin{aligned}
\dot{x}_{\theta,j} &= f(x_{\theta,j})+ G G^\top h_\theta(t,x_{\theta,j}),\;  x_{\theta,j}(T) = \xi_j, \;\; j=1,\dots,d, \\
-\dot{p}_{\theta,j} &= Df(x_{\theta,j})^\top p_{\theta,j} + D_x h_\theta(t,x_{\theta,j})^*\left(GG^\top (p_{\theta,j}+ h_\theta(t,x_{\theta,j}))\right)\\
&\quad -\alpha C^\top (y-Cx_{\theta,j}), \; \; \; p_{\theta,j}(0)=-Q_0(x_{\theta,j}(0)-x_0),\;\; j=1,\dots,d, \\
0&=\frac{1}{d}\sum_{j=1}^d\int_0^TD_\theta h_\theta(t,x_{\theta,j})^*\left( GG^\top( h_{\theta}(t,x_{\theta,j})+p_{\theta,j})\right) \; \dd t, \; \; j=1,\dots,d.
\end{aligned}
\end{equation}
For the numerical realization, we instead focus on the reduced problem and consider
\begin{align*}
 \min_{\theta \in \mathbb R^N} \frac{1}{d}\sum_{j=1}^d \left(\frac{1}{2} \| \mathcal{S}_{\theta,j}(0)-x_0 \| _{Q_0}^2 + \frac{1}{2} \int_0^T \| G^\top h_{\theta}(t,s_{\theta,j}) (t)\|^2 + \alpha \| y(t) - C\mathcal{S}_{\theta,j}(t)\|^2 \, \dd t\right),
\end{align*}
where $\mathcal{S}_{\theta,j}\colon \theta \mapsto x_{\theta_j}$ maps the vector of network parameters to the individual solutions $x_{\theta,j}$ for $j=1,\dots,d$. We utilize a gradient descent method with  Barzilai-Borwein step sizes $\gamma_\ell$ according to either one of the following rules
\begin{align*}
  \gamma_{1} = \min\left(\gamma_{\max},\frac{\langle s_{k-1},s_{k-1} \rangle }{\langle s_{k-1},y_{k-1} \rangle}\right), \quad \gamma_{2} = \min\left(\gamma_{\max},\frac{\langle s_{k-1},y_{k-1}\rangle}{\langle y_{k-1}, y_{k-1} \rangle}\right),
\end{align*}
with $s_{k-1}=\theta^{(k)}-\theta^{(k-1)}, y_{k-1}=\tilde{\theta}^{(k)} -\tilde{\theta}^{(k-1)}$ and $\tilde{\theta}^{(k)}$ given by
\begin{align*}
\tilde{\theta}^{(k)} = \frac{1}{d}\sum_{j=1}^d\int_0^TD_\theta h_{\theta^{(k)}}(t,x_{\theta^{(k)},j})^*\left( GG^\top( h_{\theta^{(k)}}(t,x_{\theta^{(k)},j})+p_{\theta^{(k)},j})\right).
\end{align*}
Here, $\gamma_{\max}$ denotes a predefined maximum stepsize.

\section{Numerical examples}

In this section, we show the numerical results obtained for the learning-based nonlinear observers and compare them with a classical (extended) Kalman filter.
For the numerical realization, we chose the network function $g_{\theta}$ to be independent of time. This is justified by the success which is achieved by autonomous networks.

All simulations were generated on an AMD Ryzen 7 1800X @ 3.68 GHz x 16,
64 GB RAM,  \matlab \;Version 9.2.0.538062 (R2017a). For the solutions of the nonlinear ODE systems, we utilize the built in  \matlab\;routine \texttt{ode15s}.

\subsection{Harmonic oscillator}

Let us consider the following undamped forced oscillator
\begin{align*}
 \ddot{x}_1(t)&=-x_1(t) + v(t), \ \  x_1(0)=x_{1,0}+\zeta_1 , \ \dot{x}_1(0) =x_{2,0}+\zeta_2 ,\\
 y(t) &= x_1(t) + w(t),
\end{align*}
where $x_1(t)$ denotes the position at time $t$ and $\zeta_1,\zeta_2,v$ and $w$ denote unknown disturbances in the initial condition, the dynamics and the observed output. Rewriting the system in first-order form yields the following system
\begin{align*}
  \dot{x}(t) &= Ax(t) + Bv(t) , \\
  y(t) & = Cx(t) + w(t),
\end{align*}
where
\begin{align*}
x(t)=
 \begin{pmatrix}
  x_1(t) \\ \dot{x_1}(t)
 \end{pmatrix},  \ A= \begin{pmatrix} 0 & 1 \\ -1 & 0 \end{pmatrix} , \ B = \begin{pmatrix} 0 \\ 1 \end{pmatrix}, \ C = \begin{pmatrix} 1 & 0 \end{pmatrix} .
\end{align*}
For generating an underlying observation $y$, we have to specify the disturbances and set
\begin{align*}
v(t)=0.1 \cos(1.2t) ,\ \  w(t)=0.1 \sin(0.5 t), \ \text{ and } \ x(0)= \begin{pmatrix} -0.1548 &     0.2969 \end{pmatrix}^\top.
\end{align*}
 For the network structure, we choose $L=2$ and $n_0,n_1,n_2=2$ such that $g_\theta$ is of the form
\begin{align*}
  g_{\theta}(z) = W_2 \left(\sigma(W_1z+b_1)+R_1z\right)
\end{align*}
with $W_1,W_2,R_1\in \mathbb R^{2\times 2}$ and $b_1\in \mathbb R^2$.

\begin{figure}[tb]
\begin{center}
  \includegraphics{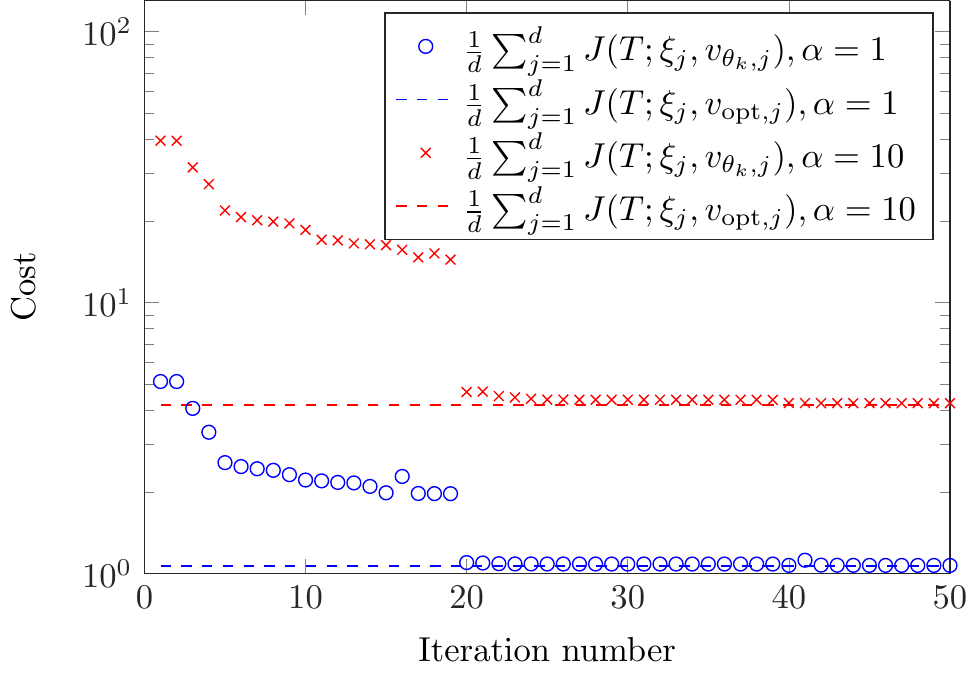}
  \caption{Behavior of the costs during the optimization for iterates $ki$ and  $d=20$ samples.}
  \label{fig:osc_cost}
  \end{center}
\end{figure}

In Figure \ref{fig:osc_cost} we show the value of the reduced cost functional $J_{\theta}$ during the optimization. The results are compared with the optimal costs that are computed by the explicit expression of the value function as in  \eqref{eq:V_linear}. Let us emphasize that the shift function $g_s(\cdot)$ was chosen such that $g_s\equiv 0$ for the iterations $k=1,\dots,20$. At iteration $k=20$, we computed a preliminary observer based on the current network function $g_\theta^{(20)}$. The associated estimate $\hat{x}_{\theta^{(20)}}$ was then used to define a shift function $g_s(t)=g_{\theta}(\hat{x}_{\theta^{(20)}})$ for the iterations $k=20,\dots,50$. Note that there is a significant decrease in the costs after the shift function $g_s$ is incorporated, indicating its importance in numerical realizations.

\begin{table}
\begin{center}
 \begin{tabular}{|l|c|c||c|c|}
 \hline
  & $J_{\mathrm{opt}},\alpha=1$ & $J_{\theta},\alpha=1$& $J_{\mathrm{opt}},\alpha=10$ & $J_{\theta},\alpha=10$\\ [0.5ex]
 \hline
 $d=1$ & 2.0606& 2.0657 & 9.0446 & 9.1329 \\
 \hline
 $d=10$ & 1.1181& 1.1244 & 4.5367 & 4.6017 \\
 \hline
 $d=20$ & 1.0672& 1.0716 & 4.1936 & 4.2517\\
 \hline
\end{tabular}
\caption{Comparison of optimal costs and network-based results.}
\label{tab:costs}
\end{center}
\end{table}

\begin{figure}[tb]
\begin{subfigure}{.5\linewidth}
\begin{center}
  \includegraphics{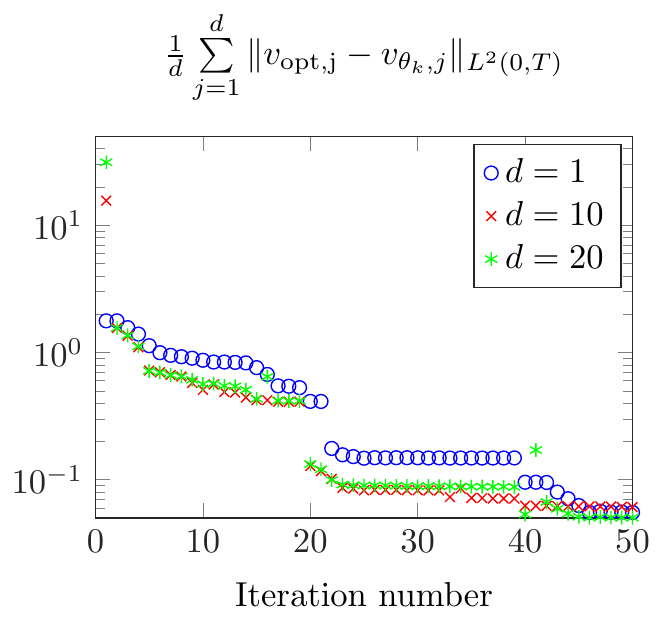}
  \caption{$\alpha=1$.}
  \label{fig:osc_l2a}
  \end{center}
\end{subfigure} \quad
\begin{subfigure}{.5\linewidth}
\begin{center}
  \includegraphics{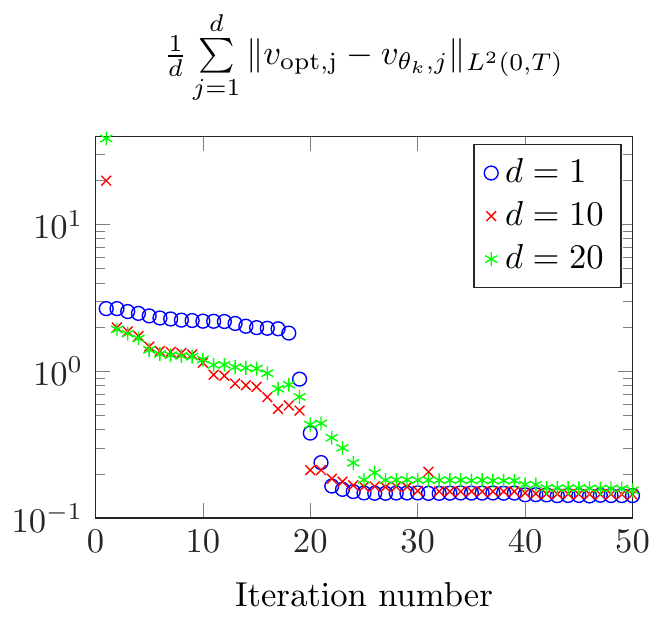}
  \caption{$\alpha=10$}
  \label{fig:osc_l2b}
  \end{center}
  \end{subfigure}
  \caption{Comparison of network-based disturbances and ``optimal disturbances''.}
  \label{fig:osc_l2}
  \quad \\
\end{figure}

From Figure \ref{fig:osc_cost} as well as Table \ref{tab:costs} it is seen that the network provides an accurate approximation of the optimal solution with respect to the value of the cost functional. This is additionally reflected in Figure \ref{fig:osc_l2} where  a comparison between network-based disturbance $v_{\theta}=G^\top h_\theta$ and ``optimal disturbance'' as provided by the Kalman filter is shown for different choices of the dimension $d$ of the sample space. Note that the optimal disturbance can be explicitly expressed (see, e.g., \cite{Willems04}) as $v_{\mathrm{opt}}(t)=\Sigma^{-1}(t)(x_{\mathrm{opt}}(t)-\hat{x}(t))$.

\begin{figure}[htb!]
\begin{subfigure}{.5\linewidth}
\begin{center}
    \includegraphics{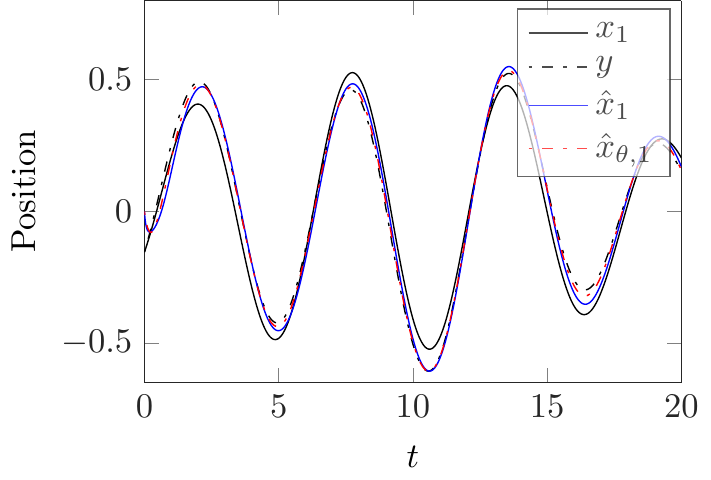}
  \end{center}
\end{subfigure} \quad
\begin{subfigure}{.5\linewidth}
\begin{center}
\includegraphics{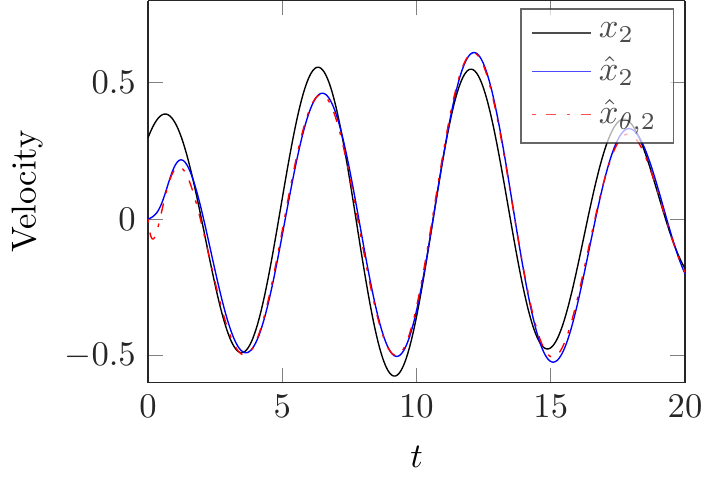}
  \end{center}
  \end{subfigure}
  \caption{Observer dynamics (network $\hat{x}_\theta$, Kalman $\hat{x}$) in comparison with ``true'' dynamics $x$ for $d=20$.}
  \label{fig:osc_obs}
\end{figure}

Finally, we compare the dynamics $\hat{x}_{\theta}$ of the network-based observer with $\hat{x}$ given by the Kalman-Bucy filter and the ``original'' dynamics $x$ that generated the observation $y$. The results are shown in Figure \ref{fig:osc_obs} and also underline the capability of $\hat{x}_{\theta}$ to estimate $x$. We observe that the Kalman filter dynamics $\hat{x}$ have the tendency to be closer to the original dynamics $x$ whereas the network-based dynamics $\hat{x}_{\theta}$ tend to be closer to the observed output $y$. A possible explanation for this behavior can be that the Kalman filter takes into account a complete probability space whereas $\hat{x}_{\theta}$ is computed on the basis of a sample space.

\subsection{Duffing oscillator}

As a nonlinear example, we consider a particular case of a general \emph{Duffing equation} for a damped and driven oscillator of the form
\begin{equation}\label{eq:duffing}
\begin{aligned}
 \ddot{x}_1(t)&+\delta \dot{x}_1(t) +\lambda x_1(t) +\beta x_1(t)^3 =v(t) \\
 y(t) &= x_1(t) + w(t).
\end{aligned}
\end{equation}
For an introduction as well as a detailed discussion of phenomena that can occur for equations of type \eqref{eq:duffing}, we refer to \cite{JorS07}. In particular, let us mention that for specific parameters configurations the dynamics of \eqref{eq:duffing} are known to exhibit chaotic behavior. Following \cite{JorS07} we therefore use the following parameters
\begin{align*}
\lambda=-1, \ \ \beta=1,\ \ \delta=0.3,\ \ v(t)=\gamma \cos(\omega t), \ \  \omega=1.2.
\end{align*}
As in the linear case, i.e., when $\delta=\beta=0$, we consider the perturbation $v$ to be a driving term for the oscillator. For the training the network parameters, we consider $\gamma =0.2$ causing a period-1 oscillation, see \cite{JorS07}. We further assume an error free measurement, i.e.,  $w \equiv 0$.  The training output $y(t)=x_1(t)$ is thus obtained on the basis of the nonlinear ODE system
\begin{align*}
 \dot{x}(t)&= Ax(t) + \begin{pmatrix} 0 \\ - x_1(t)^3 \end{pmatrix} + \begin{pmatrix} 0\\ 0.2 \end{pmatrix} v(t), \ \ x(0) = \begin{pmatrix} 0.0646  &   -0.1465 \end{pmatrix}^\top.
\end{align*}

\begin{figure}[htb]
\begin{subfigure}{.5\linewidth}
\begin{center}
    \includegraphics{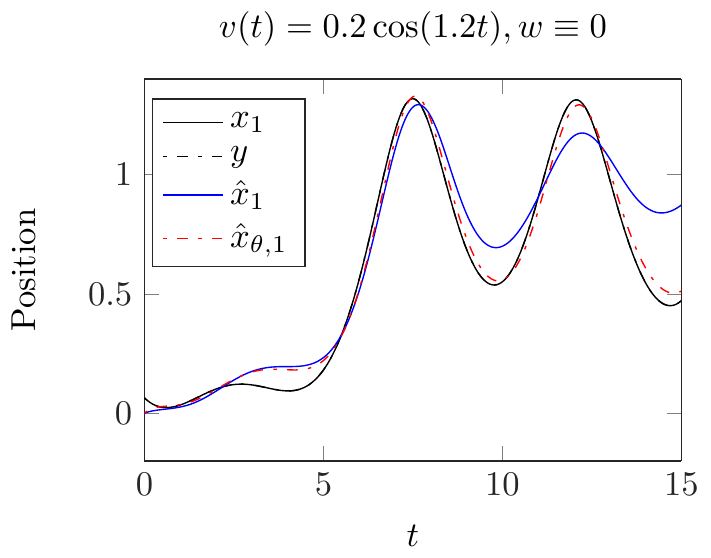}
  \end{center}
\end{subfigure} \quad
\begin{subfigure}{.5\linewidth}
\begin{center}
  \includegraphics{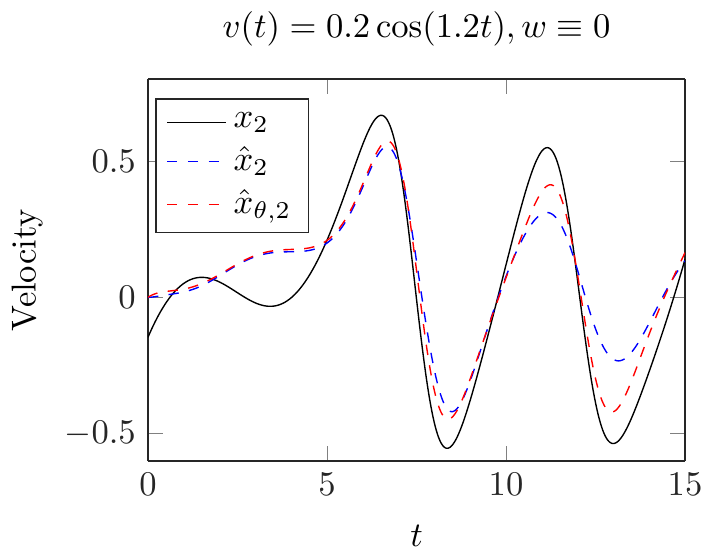}
  \end{center}
  \end{subfigure}
  \caption{Observer dynamics (network $\hat{x}_{\theta}$, extended Kalman $\hat{x}$) in comparison with training data $x$ for $\alpha=1$ and $d=5$.}
    \label{fig:duff_train}
    \quad \\
\end{figure}

In Figure \ref{fig:duff_train}, we show the results of a network-based observer on the basis of the learning problem \eqref{eq:Ptheta} for $d=5$ samples $\xi_j$ (randomly chosen) and $\alpha=1$ on the time interval $[0,15]$. We also depict  the  estimate associated with the extended Kalman filter which is based on a recursive update of the observer gain $\Sigma_{\hat{x}(t)} C^\top$ for the linearized dynamics along the current state estimate $\hat{x}$. For $x_1(\cdot)=y(\cdot)$, the estimate of the network-based observer is in very good agreement with the original dynamics, while for $x_2(\cdot)$ the results are slightly worse both for the extended Kalman as well as the network-based estimate which is nevertheless performing better. For larger values of $\alpha$, we observed an (expected) improvement of both approximations, with the network-based observer adjusting better than the extended Kalman filter.

\begin{figure}[tb]
\begin{subfigure}{.5\linewidth}
\begin{center}
  \includegraphics{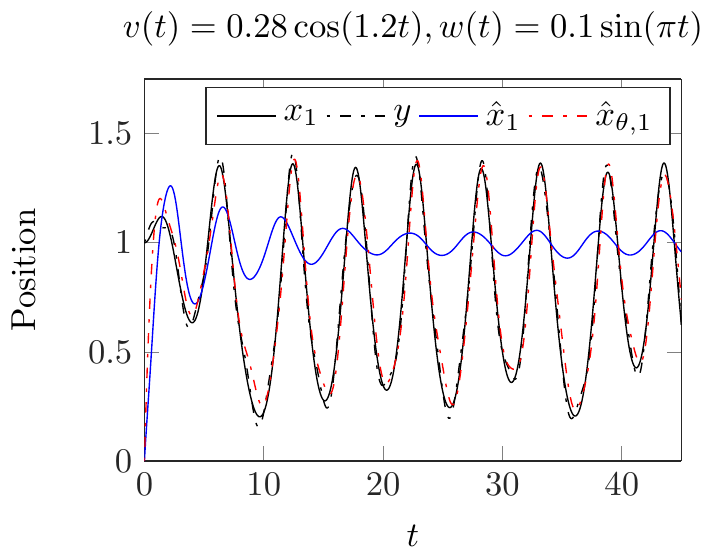}
  \end{center}
\end{subfigure} \quad
\begin{subfigure}{.5\linewidth}
\begin{center}
  \includegraphics{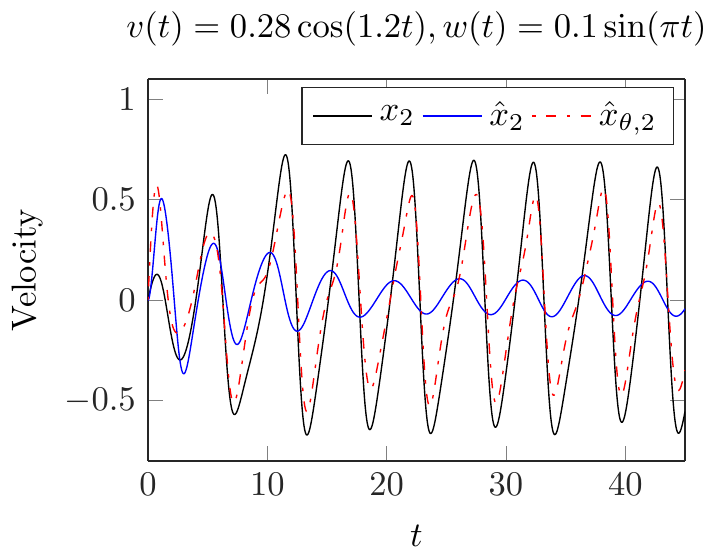}
  \end{center} 
  \end{subfigure}
  \caption{Observer dynamics (network $\hat{x}_{\theta}$, extended Kalman $\hat{x}$) in comparison with test data $x$.}
  \label{fig:duff_test}
\end{figure}

With regard to the discussion in Section \ref{subsec:y-indp}, in Figure \ref{fig:duff_test} we also show the results of both observers for a different driving term $v(\cdot)$, causing a period-2 oscillation, as well as an additional measurement error $w(\cdot)=0.1\sin(\pi \cdot)$. We stress that the network parameters $\theta$ remained to be those obtained in the training step and were not adjusted to the new data set.

We emphasize that the observation $y$ not only differs since different $v$ and $w$ were used, but also the initial condition $x_0$, as well as the time horizon $T$ were changed. As is apparent from Figure \ref{fig:duff_test}, the network-based observer is still reproducing  the qualitative behavior of the new dynamics $y=x_1+w$ and $\dot{x}_1$. On the other hand, the state estimate utilizing the extended Kalman filter shows significant deviations in this case.


\begin{remark}
Let us point out that the very formulation of the Mortensen observer in \eqref{eq:P_T}  involves solving   the underlying dynamical system backwards in time. The necessity for  backward solves also appears in the optimality system \eqref{eq:network_1st_opt}.  It may lead to difficulties for systems with a special structure, as for instance the \emph{Van der Pol oscillator} whose dynamics is of the form
\begin{align*}
 \ddot{x}_1-(1-x_1^2)\dot{x}_1+x_1 =0.
\end{align*}
This system exhibits a repulsive limit cycle when it is considered
backwards in time  see, e.g., \cite{MarRF15}. This  causes numerical instabilities for a large class of ODE solvers, see \cite{Haf19}. In particular, once the trajectory is outside the limit cycle it rapidly escapes to infinity, and if numerical error happens along the trajectory inside the limit cycle the numerical solution may tend to the origin prematurely. For such systems, the Mortensen observer may not be the method of choice for state reconstruction. 
\end{remark}

\bibliographystyle{siam}
\bibliography{references}

\begin{thebibliography}{10}

\bibitem{Ad12}
{\sc D.~M. Adhyaru}, {\em State observer design of nonlinear systems using
  neural networks}, Applied soft computing, 12 (2012), pp.~2530--2537.

\bibitem{BucJ68}
{\sc R.~S. Bucy and P.~D. Joseph}, {\em Filtering for stochastic processes with
  applications to guidance}, Interscience Tracts in Pure and Applied
  Mathematics, No. 23, Interscience Publishers John Wiley \& Sons., Inc., New
  York-London-Sydey, 1968.

\bibitem{Fleming97}
{\sc W.~H. Fleming}, {\em Deterministic nonlinear filtering}, Annali della
  Scuola Normale Superiore di Pisa - Classe di Scienze, Ser. 4, 25 (1997),
  pp.~435--454.

\bibitem{FS06}
{\sc W.~H. Fleming and H.~M. Soner}, {\em Controlled {M}arkov processes and
  viscosity solutions}, vol.~25 of Stochastic Modelling and Applied
  Probability, Springer, New York, second~ed., 2006.

\bibitem{GKNV19}
{\sc R.~Gribonval, G.~Kutyniok, M.~Nielsen, and F.~Voigtlaender}, {\em
  Approximation spaces of deep neural networks}, 2019.
\newblock available from \url{https://arxiv.org/abs/1905.01208}.

\bibitem{GPEB19}
{\sc P.~Grohs, D.~Perekrestenko, D.~Elbr\"achter, and H.~B\"olcskei}, {\em Deep
  neural network approximation theory}, tech. rep., 2019.
\newblock available from \url{https://arxiv.org/abs/1901.02220}.

\bibitem{Haf19}
{\sc S.~Hafstein}, {\em Numerical Analysis Project in ODEs for Undergraduate
  Students}, 06 2019, pp.~421--434.

\bibitem{HeZRS16}
{\sc K.~{He}, X.~{Zhang}, S.~{Ren}, and J.~{Sun}}, {\em Deep residual learning
  for image recognition}, in 2016 IEEE Conference on Computer Vision and
  Pattern Recognition (CVPR), June 2016, pp.~770--778.

\bibitem{Jaz70}
{\sc A.~H. Jazwinski}, {\em Stochastic processes and filtering theory}, Courier
  Corporation, 1970.

\bibitem{JorS07}
{\sc D.~Jordan and P.~Smith}, {\em Nonlinear ordinary differential equations:
  an introduction for scientists and engineers}, vol.~10, Oxford University
  Press on Demand, 2007.

\bibitem{KalB61}
{\sc R.~E. Kalman and R.~S. Bucy}, {\em New results in linear filtering and
  prediction theory}, Transactions of the ASME. Series D. Journal of Basic
  Engineering, 83 (1961), pp.~95--108.

\bibitem{KK85}
{\sc H.-W. Knobloch and H.~Kwakernaak}, {\em Lineare {K}ontrolltheorie},
  Springer-Verlag, Berlin, 1985.

\bibitem{Krener04}
{\sc A.~J. Krener}, {\em The convergence of the minimum energy estimator}, in
  New Trends in Nonlinear Dynamics and Control and their Applications, W.~Kang,
  C.~Borges, and M.~Xiao, eds., Berlin, Heidelberg, 2003, Springer Berlin
  Heidelberg, pp.~187--208.

\bibitem{KW20}
{\sc K.~Kunisch and D.~Walter}, {\em Semiglobal optimal feedback stabilization
  of autonomous systems via deep neural network approximation}, tech. rep.,
  2020.
\newblock available from \url{https://arxiv.org/abs/2002.08625}.

\bibitem{MarRF15}
{\sc A.~Marton, A.~Ribeiro, and A.~Fioravanti}, {\em Comparison between {SOS}
  and {(S)DSOS L}yapunov functions for nonlinear systems}, 12 2015.

\bibitem{Moireau18}
{\sc {Moireau, P.}}, {\em A discrete-time optimal filtering approach for
  non-linear systems as a stable discretization of the {M}ortensen observer},
  ESAIM: Control, Optimisation and Calculus of Variations, 24 (2018),
  pp.~1815--1847.

\bibitem{Mortensen68}
{\sc R.~E. Mortensen}, {\em Maximum-likelihood recursive nonlinear filtering},
  Journal of Optimization Theory and Applications, 2 (1968), pp.~386--394.

\bibitem{Son98}
{\sc E.~Sontag}, {\em Mathematical Control Theory: Deterministic Finite
  Dimensional Systems}, vol.~6, Springer Verlag, 1998.

\bibitem{Willems04}
{\sc J.~Willems}, {\em Deterministic least squares filtering}, Journal of
  Econometrics, 118 (2004), pp.~341 -- 373.
\newblock Contributions to econometrics, time series analysis, and systems
  identification: a Festschrift in honor of Manfred Deistler.

\end{thebibliography}

\subsection*{Acknowledgement}

This work was  supported in  part by the ERC advanced grant 668998 (OCLOC) under
the EU's H2020 research program. We would like to thank Daniel Walter (RICAM, Linz) for several fruitful discussions on the topic.

\end{document}